\documentclass[11pt]{amsart}

\usepackage[USenglish]{babel}

\usepackage{amsmath,amsthm,amssymb,amscd}
\usepackage{booktabs}
\usepackage{url}

\usepackage{MnSymbol}

\usepackage{tikz}
\usetikzlibrary{arrows,shapes.geometric,positioning,decorations.markings}

\usepackage[mathscr]{eucal}
\usepackage[normalem]{ulem}
\usepackage{latexsym,youngtab}
\usepackage{multirow}
\usepackage{epsfig}
\usepackage{parskip}
\usepackage[all,cmtip]{xy}

\usepackage[margin=2cm]{geometry}
\usepackage{color}

\newtheoremstyle{custom}
  {3pt}
  {3pt}
  {\slshape}
  {}
  {\bfseries}
  {.}
  { }
   {}
\theoremstyle{custom}
\newtheorem{theorem}{Theorem}[section]
\newtheorem{proposition}[theorem]{Proposition}

\newtheorem{proposition/definition}[theorem]{Proposition/Definition}

\theoremstyle{definition}
\newtheorem{definition}[theorem]{Definition}

\newtheorem{example}[theorem]{Example}

\theoremstyle{remark}
\newtheorem{remark}[theorem]{Remark}

\newtheoremstyle{exercise}
  {3pt}
  {6pt}
  {}
  {}
  {\bfseries}
  {:}
  { }
   {}
\theoremstyle{exercise}
\newtheorem{exercise}[theorem]{Exercise}
\newtheoremstyle{exercises}
  {3pt}
  {6pt}
  {}
  {}
  {\bfseries}
  {:}
  {\newline}
   {}

\theoremstyle{exercise}
\newtheorem{exercises}[theorem]{Exercises}

\input epsf
\def\boxit#1{\vbox{\hrule height1pt\hbox{\vrule width1pt\kern3pt
  \vbox{\kern3pt#1\kern3pt}\kern3pt\vrule width1pt}\hrule height1pt}}

\def\trank{\text{rank}}

\def\BC{\mathbb C}

\def\BP{\mathbb P}
\def\pp#1{\mathbb P^{#1}}

\def\fgl{\mathfrak g\mathfrak l}

\def\pp#1{{\mathbb P}^{#1}}
\def\tdim{{\rm dim}}

\def\ww{\wedge}

\def\inv{{}^{-1}}

\def\cA{{\mathcal A}}

\def\cO{{\mathcal O}}

\def\11{\mathbf 1}

\def\fg{{\mathfrak g}}

\def\a{\alpha}

\def\b{\beta}
\def\g{\gamma}
\def\s{\sigma}

\def\ot{{\mathord{ \otimes } }}
\def\op{{\mathord{\,\oplus }\,}}

\def\ra{{\mathord{\;\rightarrow\;}}}

\def\La#1{\Lambda^{#1}}

\def\frak{\mathfrak}
\def\fgl{\frak g\frak l}

\def\op{\oplus}

\def\ep{\epsilon}
\def\op{\oplus}

\def\ul{\underline}

\def\s{\sigma}
\def\t{\tau}

\def\a{\alpha}
\def\b{\beta}

\def\g{\gamma}

\def\FS{\mathfrak  S}

\def\ol{\overline}

\def\BP{\mathbb  P}
\def\BC{\mathbb  C}

\def\pp#1{\mathbb  P^{#1}}

\def\ep{\epsilon}

\def\ci{\mathcal  I}

\def\fg{\mathfrak  g}

\def\hd{, \dotsc ,}

\def\inv{{}^{-1}}

\def\La#1{\Lambda^{#1}}

\def\pp#1{\mathbb  P^{#1}}

\def\ur{\underline{\mathbf{R}}}
\def\uR{\underline{\mathbf{R}}}

\def\ra{\rightarrow}

\def\tend{\operatorname{End}}

\def\tdim{\operatorname{dim}}

\def\tlim{\lim}
\def\tmod{\operatorname{mod}}

\def\tmin{\operatorname{min}}

\def\trank{\operatorname{rank}}

\def\ww{\wedge}

\def\be{\begin{equation}}
\def\ene{\end{equation}}
\def\aaa{{\mathbf{a}}}

\DeclareMathOperator{\tlog}{log}

\def\trank{\mathbf{R}}

\newcommand{\Id}{\operatorname{Id}}

\def\Mn{M_{\langle \nnn \rangle}}\def\Mone{M_{\langle 1\rangle}}

\def\Mn{M_{\langle \nnn \rangle}}\def\Mone{M_{\langle 1\rangle}}

\def\trank{{\mathrm {rank}}}

\def\aaa{\mathbf{a}}

\def\nnn{\mathbf{n}}

\def\Mn{M_{\langle \nnn \rangle}}\def\Mone{M_{\langle
1\rangle}}

\def\lam{\lambda}

\def\trank{{\mathrm {rank}}}

\def\aaa{{\bold a}}

\def\nnn{\bold n}

\renewcommand{\a}{\alpha}
\renewcommand{\b}{\beta}
\renewcommand{\g}{\gamma}

\renewcommand{\BC}{\mathbb{C}}

\newcommand{\bfR}{\mathbf{R}}
  
\renewcommand{\bar}[1]{\overline{#1}}
\renewcommand{\hat}[1]{\widehat{#1}}

 \renewcommand{\tilde}{\widetilde}

\def\om{\omega}

   \def\Amat{X}
        \def\Bmat{Y}
        \def\Cmat{Z}
     \newcommand{\acta}{\circ_{\scriptscriptstyle A}}
    \newcommand{\actb}{\circ_{\scriptscriptstyle B}}
    \newcommand{\actc}{\circ_{\scriptscriptstyle C}}
    
      \newcommand{\alg}[1]{\cA_{111}^{#1}} 
 
\begin{document}

\author{J. M. Landsberg}

\address{Department of Mathematics, Texas A\&M University, College Station, TX 77843-3368, USA}
 
\email{jml@math.tamu.edu}

\title[Secant varieties and matrix multiplication]{Secant varieties and the complexity of matrix multiplication}

\thanks{Landsberg   supported by NSF grants  AF-1814254 and  and AF-2203618}



\begin{abstract}  
  This is a survey primarily about determining the border rank of tensors, especially those relevant for the
  study of the complexity of matrix multiplication. 
  This is a subject that on the one hand is of great significance in theoretical computer
  science, and on the other hand touches on many beautiful topics in algebraic geometry such as 
 classical and recent results on equations for secant varieties (e.g., via vector bundle 
 and representation-theoretic methods) and the geometry and deformation theory of
 zero dimensional schemes.
\end{abstract}

\maketitle

\section{Introduction}
This is a survey of uses of secant varieties in the study of the complexity of  matrix multiplication, one of many
areas in which    Giorgio Ottaviani has made significant contributions.
I pay special attention to   the use of deformation theory 
because at this writing, deformation theory provides the most promising path  to overcoming lower bound barriers. For an introduction
to more general uses of algebraic geometry  in algebraic complexity theory see \cite{alexsurvey}.
I begin by reviewing  some classical results.

\subsection{Symplectic bundles on the plane, secant varieties, and 
L\"uroth quartics revisited  \cite{MR2554725}}

In the 1860's, Darboux   studied degree $n$ curves in $\pp 2$ that pass through all the $\binom {n+1}2$ vertices  of a complete $(n+1)$-gon
in $\pp 2$ (i.e., the union of $n+1$ lines in $\pp 2$ with no points of triple intersection).   
In 1869  
L\"uroth studied the $n=4$ case.  A na\"\i ve dimension count   indicates
that  all quartics  should   pass through   the $10$ vertices  of some complete pentagon
  but L\"uroth proved it is actually a    codimension one condition.

 In  
1902
Dixon \cite{dixon}   proved all degree $n$ curves in $\pp 2$ arise as a $n\times n$ symmetric determinant
(also see \cite{MR1501168} for the general determinantal case).   

In 1977
Barth \cite{MR0460330} studied the moduli space of stable (symplectic) vector bundles on $\pp 2$.
In particular he showed that
the curve of jumping lines of a rank $2$ stable bundle
on
$\BP^2$
 with Chern classes $(c_1, c_2) = (0, 4)$ is a L\"uroth quartic.   Barth
also gave   a new proof of L\"uroth's theorem via vector bundles.

In \cite{MR2554725} Giorgio Ottaviani  explains these  results via the defectivity of  {\it secant varieties} of $Seg(\pp 2\times v_2(\pp{n-1}))$, where $Seg(\pp 2\times v_2(\pp{n-1}))\subset \BP (\BC^3\ot S^2\BC^n)$ is the set of points
$[x\ot z^2]$, where $[x]\in \pp 2$ and $z\in\pp{n-1}$.
The proof    
uses the bounded derived category version of   Beilinson's monad Theorem \cite{MR509387}, see
\cite{MR1113077} for an excellent introduction.

\subsection{Secant varieties}  
Throughout this paper $V,A,B,C$ denote finite dimensional complex vector spaces.
Let
$X\subset \BP V$ be a  projective variety, 
  Define its {\it $r$-th secant variety},  or variety of secant $\pp{r-1}$'s,  to be
$$\s_r(X):=\ol{\bigcup_{x_1\hd x_r\in X} \langle x_1\hd x_r\rangle }.
$$
Here, for a set or subscheme $Z\subset \BP V$, $\langle Z\rangle\subset \BP V$ denotes its linear span,
and the overline denotes Zariski closure.

In this article I will be particularly interested in the case $X=Seg(\BP A\times \BP B\times\BP C)
\subset \BP (A\ot B\ot C)$,  the variety of rank one ($3$-way) tensors. Given $T\in A\ot B\ot C$,
define the {\it border rank} of $T$, $\ur(T)$ to be the smallest $r$ such that $[T]\in \s_r(X)$.

Secant varieties have a long history in algebraic geometry dating back to the 1800's. In the 20th
century they were used by
J. Alexander and A. Hirschowitz \cite{AH} to solve the polynomial interpolation problem,
and by F. Zak \cite{Zak} to solve a linearized version of R. Hartshorne's famous conjecture
on complete intersections, called {\it Hartshorne's conjecture on linear normality}.
L. Manivel and I used them to study the geometry of the exceptional groups and their
homogeneous varieties,  and even to obtain  a new proof of the Killing-Cartan
classification of complex simple Lie algebras and prove geometric consequences
of conjectured categorical generalizations of Lie algebras by Deligne and Vogel, see    \cite{MR2090671} for a survey.
In this article, I  discuss their use in the context of {\it algebraic complexity theory},
more specifically, in proving lower  and upper  bounds on the {\it complexity of matrix multiplication}.

\subsection{Matrix multiplication}
In   1968, V. Strassen \cite{Strassen493}  discovered the usual way we multiply $\nnn\times\nnn$-matrices,  which
uses $\cO(\nnn^3)$ arithmetic operations, is not optimal. After much work, it was generally conjectured
that one can in fact multiply matrices using $\cO(\nnn^{2+\ep})$ arithmetic operations for any $\ep>0$.
To fix ideas, 
define the {\it exponent of matrix multiplication} $\om$ to be the infimum over all $\t$ such that
$\nnn\times\nnn$ matrices may be multiplied using $\cO(\nnn^\t)$ arithmetic operations,
so the conjecture is that $\om=2$. 
The {\it matrix multiplication tensor}  $\Mn: \BC^{\nnn^2}\times \BC^{\nnn^2}\ra \BC^{\nnn^2}$
  executes the bilinear map of multiplying two matrices.
Fortunately for algebraic geometry,
Bini \cite{MR605920} showed $ \ur(\Mn) =\cO(\nnn^\om)$ so we may study the exponent via secant
varieties of Segre varieties.

Thus one way to prove complexity {\it lower} bounds   for matrix multiplication would be  to prove lower bounds on 
the border rank of $\Mn$.   I will give a history of such lower bounds.
Perhaps more surprising, is that one way of  showing 
{\it upper} bounds  for the complexity of matrix multiplication would be to prove
the border rank of certain auxiliary tensors is small, as I discuss in \S\ref{upperbnds}.

\subsection{Dimensions of secant varieties}\label{secdimsect}
One expects $\tdim \s_r(X)=\tmin\{ r\tdim X+r-1, \tdim \BP V\}$, because one can pick $r$ points on $X$, and
then a point on the $\pp{r-1}$ spanned by them.   This always gives an upper bound on the dimension. 

Strassen \cite{Strassen505}, motivated by  the complexity of matrix multiplication,  showed that this  expectation fails for  $X=Seg(\pp 2\times \pp {n-1}\times \pp{n-1})
\subset \BP (\BC^3\ot \BC^n\ot \BC^n)$, $n$ odd, $r=\frac{3n-1}2$.
   
Previously, 
E. Toeplitz, in  1877 \cite{MR1509924}, had already
shown it    fails for   $X=Seg(\pp 2\times v_2(\pp{3}))\subset \BP (\BC^3\ot S^2\BC^4)$, $r=5$.

In  2007 Ottaviani  \cite{MR2554725}  showed that  more generally
the expectation  fails for  $X=Seg(\pp 2\times v_2(\pp{n-1}))\subset \BP (\BC^3\ot S^2\BC^n)$ with  $n$ even
$r=\frac{3n}2-1$, and that  this failure implies L\"uroth's theorem.   
In the same paper he also partially recovers Barth's moduli results.

\subsection{Acknowledgements} I thank the organizers of GO60 for putting together a wonderful conference under
difficult conditions. I also thank J. Jelisiejew and M. Micha{\l}ek for permission to include section  \ref{endclosedpf}, which arose
out of conversations with them. Most of all I thank Giorgio Ottaviani for our collaborations together that I hope continue
into the future.
   
\section{Koszul flattenings and variants} 
\subsection{Idea of Proofs of results in \S\ref{secdimsect}}  
To prove the na\"\i ve dimension count for $\tdim \s_r(X)$ is wrong (e.g.,   in the case of L\"uroth's theorem
that $\s_r(X)\neq \BP V$),
one can
show that  the  ideal of $\s_r(X)$ is non-empty by exhibiting an explicit polynomial in the ideal.

Strassen did this and his result was  revisited by Ottaviani:
Consider  $Seg(\BP A\times \BP B\times \BP C)\subset \BP (A\ot B\ot C)$,
$\tdim A=3$, $\tdim B=\tdim C=m$.
Let $\{a_i\},\{ b_j\},\{ c_k\}$ be bases of $A,B,C$. Given $T=\sum T^{ijk}a_i\ot b_j\ot c_k \in A\ot B\ot C$, 
consider the linear map
\begin{align*}
T^{\ww 1}_A: A\ot B^*&\ra \La 2 A\ot C
\\
 a\ot \beta&\mapsto \sum_{i,j,k} T^{ijk}\b(b_j) a\ww a_i\ot c_k
\end{align*}

Exercise: If  $[T]\in Seg(\BP A\times \BP B\times \BP C)$, then  $\trank (T^{\ww 1}_A)=2$, and thus,
by linearity,  
  if   $\trank(T^{\ww 1}_A)> 2R$, then  $[T]\not\in \s_R (Seg(\BP A\times \BP B\times \BP C))$.

Ottaviani states in a remark that these minors are a reformulation of {\it Strassen's equations}
(however see Remark \ref{p=1stronger} below),
which, for tensors $T\in A\ot B\ot C=\BC^\aaa\ot \BC^m\ot \BC^m$ such that there exists     $\a \in A^*$
with   $\trank( T(\a))=m$,
are naturally expressed as follows:
consider  $T(A^*) \subset B\ot C$,   and for   $\a \in A^*$  with  $\trank (T(\a))=m$,  consider 
the linear isomorphism   $T(\a): B^*\ra C$.
Then  $T(A^*)T(\a)^{\inv}\subset \tend(C)$  is a space of endomorphisms.
If 
$T=\sum_{j=1}^m e_j\ot b_j\ot c_j$ for some $e_j\in A$, then one obtains a space of diagonal
matrices. In particular, the matrices commute. Since the property of commuting is closed, if
$[T]\in \s_m(Seg(\BP A\times \BP B\times \BP C))$, then one still obtains a space of
 commuting endomorphisms.   
Moreover, the rank of the commutator (a measure of the failure
of commutivity) may be computed from the  rank of  $T^{\ww 1}_A$.  
Note that in both cases one restricts to a three dimensional subspace of $A$.

To see Strassen's equations as polynomials, for $X\in B\ot C$, let $ X^{\ww m-1}\in \La{m-1}B\ot \La{m-1}C\simeq B^*\ot C^*$
denote the adjucate (cofactor matrix), and recall that $X\inv$ is essentially the adjugate times the determinant. 
Then Strassen's equations for $T$ to have border rank (at most)  $m$ \cite{Strassen505} become, for all
 $X,Y,Z\in T(A^*)\subset B\ot C$,
\be\label{strassensubspace}
X Y^{\ww m-1} Z-Z Y^{\ww m-1}X=0.
\ene
These are equations of degree $m+1$.

Using a refinement of these equations that
takes into account the rank of the commutator, Strassen proved $\ur(\Mn)\geq \frac 32 \nnn^2$, the first non-classical lower bound
on the border rank of the matrix multiplication tensor.

Call a tensor $T$ which  satisfies the genericity condition   that there exists  $\a \in A^*$ with $\trank (T(\a))=m$,
 {\it $1_A$-generic}.

When $T$ is $1_A$-generic, taking $Y$ of full rank and changing bases such that it is the identity element,
the equations require the space to be abelian. If $T(A^*)$ is of bounded rank $m-1$, for each $X,Y,Z$,
the set of $m^2$ equations reduces to a single equation. If $T(A^*)$ is of bounded rank $m-2$, then
the equations become vacuous. 

 Lemma 3.1 of \cite{MR2554725} states that if $T$ is $1_A$-generic, then the condition on 
$\trank(T^{\ww 1}_A)$ is indeed a reformulation of Strassen's equations.

 \begin{remark} \label{p=1stronger}
 Recently, with my student Arpan Pal and Joachim Jelisiejew \cite{https://doi.org/10.48550/arxiv.2205.05713}, we proved that
   if $T$ is  not $1_A$-generic,  then   the condition on  $\trank(T^{\ww 1}_A)$ is a  stronger condition than the $A$-Strassen equations. 
 \end{remark}

\subsection{Generalizations: Young flattenings  \cite{MR3081636,MR3376667} }

\subsubsection{Koszul flattenings}
Consider  $Seg(\BP A\times \BP B\times \BP C)\subset \BP (A\ot B\ot C)$,
let $\tdim A=2p+1$, $\tdim B=\tdim C=m$. (If $\tdim A>2p+1$, restrict to a general $2p+1$ dimensional
subspace.)
  
Given $T=\sum T^{ijk}a_i\ot b_j\ot c_k \in A\ot B\ot C$, 
consider the linear map
\begin{align*}
T^{\ww p}_A: \La p A\ot B^*&\ra \La {p+1} A\ot C
\\
 a_{i_1}\ww\cdots \ww a_{i_p}\ot \beta&\mapsto \sum_{i,j,k} T^{ijk}\b(b_j) a_{i_1}\ww\cdots \ww a_{i_p}\ww a_i\ot c_k
\end{align*}

Exercise: If  $[T]\in Seg(\BP A\times \BP B\times \BP C)$, then  $\trank (T^{\ww p}_A)=\binom{2p}{p-1}$.    
Thus if   $\trank(T^{\ww p}_A)>  \binom{2p}{p-1}R$, then  $[T]\not\in \s_R (Seg(\BP A\times \BP B\times \BP C))$. 
Call these equations the {\it $p$-Koszul flattenings}.

When Ottaviani and I  found the $p$-Koszul flattenings, we were sure we would get a new lower bound for matrix
multiplication. Our first attempts were discouraging, we were attempting to take $2p+1=\nnn^2$
or nearly so. It turns out that our initial attempts were too greedy, as such values do not give
good lower bounds. Only months later, we finally tried taking $p=  \nnn-1 $  
which enabled us  to prove the first new  lower bounds for border rank of matrix multiplication since 1983:

\begin{theorem} \cite{MR3376667}  $\ur(\Mn)\geq 2\nnn^2-\nnn$.
\end{theorem}

It is worth noting that the absolute limit of this method, and any determinantal equations that we found, was
  $2\nnn^2-1$, i.e., for tensors in $\BC^m\ot \BC^m\ot \BC^m$,    $2m-1$.

\subsubsection{Young flattenings}
We found  the $p$-Koszul flattenings as part of a general program to systematically
find equations for secant varieties, especially equations of secant varieties of homogeneous varieties,
which we call {\it Young flattenings}.
Giorgio likes to think of these in terms of degeracy loci of maps between vector bundles, and I prefer
a more representation-theoretic perspective. The basic idea is for $X\subset \BP V$,
to find an inclusion of $V$ into a space of matrices. Then if the matrices
associated to points of $X$ have rank at most $q$, the size $qr+1$ minors restricted
to $V$ give equations for $\s_r(X)$.  

\subsubsection*{Vector bundle perspective}
   Let $E\ra X$ be a vector bundle of rank $e$,
write $L=\cO_X(1)$ so $V=H^0(X,L)^*=H^0(L)^*$. Let $v\in V$ and consider the linear map
$$
A^E_v: H^0(E)\ra H^0(E^*\ot L)^*
$$
induced by the multiplication map $H^0(E)\ot H^0(L)^*\ra H^0(E^*\ot L)^*$. Then, assuming
all spaces are sufficiently large,  the size
$(re+1)$ minors of $A^E_v$ give equations for $\s_r(X)$. Here we have an
inclusion $V=H^0(L)^*\subset H^0(E)^*\ot H^0(E^*\ot L)^*$.

\subsubsection*{Representation theory perspective}
Let $X=G/P\subset \BP V_{\lam}$ where $V_{\lam}$ is an irreducible module
for the reductive group $G$ of highest weight $\lam$ and $X$ is the orbit of a highest weight line, i.e., the minimal
$G$-orbit in $\BP V_{\lam}$. Look for $G$-module inclusions $V_{\lam}\subset V_{\mu}\ot V_{\nu}$,
so in coordinates one realizes $V_{\lam}$ as a space of matrices. Then for  $x\in X$ if the associated
matrix has rank $k$, the size $rk+1$ minors of $V_{\mu}\ot V_{\nu}$ restricted to $V_{\lam}$
give equations in the ideal of $\s_r(X)$.

We spent some time trying to find more powerful Young flattenings. 
At first we just thought we were not being clever enough in our search for determinantal equations,
but then we came to suspect that there was some limit to the method.

\section{Beyond Koszul flattenings: steps forward and barriers to future progress}
\subsection{The cactus barrier}
Around the same time, in both algebraic complexity theory \cite{MR3761737}
and algebraic geometry \cite{MR2996880,MR3121848} (see \cite[Chap. 10]{MR3729273} for an overview),
it was proven that determinantal methods are subject to an absolute barrier that is at most $6m$
for tensors in $\BC^m\ot \BC^m\ot \BC^m$.

To explain the barrier from a geometric perspective, rewrite the definition of the secant variety as
$$\s_r(X):=\overline{ \bigcup  \{ \langle R\rangle\mid {\rm length}(R)=r,\ 
 R \subset X,\ R:{\rm smoothable} \} }.
$$
Here $R\subset X$ denotes a zero dimensional scheme and the union is taken over all zero dimensional
schemes satisfying the conditions in braces.
Define the {\it cactus variety} \cite{MR3121848}:
$$\kappa_r(X):=\overline{ \bigcup  \{ \langle R\rangle\mid {\rm length}(R)=r,\ 
 R \subset X\} }.
$$

It turns out that $\kappa_{6m}(Seg(\pp{m-1}\times\pp{m-1}\times\pp{m-1}))=\BP(\BC^m\ot \BC^m\ot \BC^m)$, 
compared with $\s_r(Seg(\pp{m-1}\times\pp{m-1}\times\pp{m-1}))$ which does
not fill $\BP(\BC^m\ot \BC^m\ot \BC^m)$
until $r\sim \frac{m^2}3$.

The barrier results from the fact that {\it determinantal equations for $\s_r(X)$ are
also equations for $\kappa_r(X)$}!

When I learned this, I became very discouraged.

\subsection{A Phyrric victory}
With M. Micha{\l}ek, we were able to push things a little further for tensors with symmetry.
Given
$T\in A\ot B\ot C$,
$\ur(T)\leq r$ if and only if there
exists a  curve  $E_t\subset G(r,A\ot B\ot C)$
such that 
\begin{itemize}

\item For $t\neq 0$, $E_t$ is spanned by $r$ rank one elements. 

\item $T\in E_0$.
\end{itemize}

Let $G_T:=\{ g\in GL(A)\times GL(B)\times GL(C) \mid g\cdot T=T\}$ denote the symmetry group of $T$.
Then if we have such a curve $E_t$, then for all  $g\in G_T$, $gE_t$ also
gives a border rank decomposition.
Thus one
  can insist on normalized curves, those with $E_0$ Borel fixed
  for a Borel subgroup of $G_T$  \cite{MR3633766}. Then one can apply a border rank version
  of the classical substitution method (see, e.g., \cite{MR3025382}) to reduce the problem
  to bounding the border rank of a smaller tensor.
Applying this to the matrix multiplication tensor, we proved:

\begin{theorem}
\cite{MR3842382}  $\ur(\Mn)\geq 2 \nnn^2-\tlog_2\nnn-1$.
\end{theorem}

Recall that the limit of lower bounds one can prove with  Young flattening  is $2m-1$.
We   wrote down an explicit sequence of tensors $T_m\in \BC^m\ot \BC^m\ot \BC^m$
with a one-dimensional symmetry group
and proved:

\begin{theorem}\cite{landsberg2019finding} $\ur(T_m)\geq (2.02)m$.
\end{theorem}

 After that, I did not see any path to further lower bounds.

\subsection{A vast generalization: border apolarity}
W. Buczy{\'n}ska and J. Buczy{\'n}ski \cite{MR4332674} had the following idea: Consider not just a 
curve in the Grassmannian obtained by taking the spans of $r$  moving  points $\{T_{j,\ep}\}$, where
$T=\tlim_{\ep\ra 0} \sum_{j=1}^r T_{j,\ep}$,  
 but the curve of {\bf ideals}
$I_\ep\in Sym(A^*\op B^*\op C^*)$ that the points  give  rise to: 
let 
$I_\ep$ be the  ideal of polynomials vanishing on 
 $[T_{1,\ep}]\cup \cdots \cup [T_{r,\ep}]\subset \BP A\times \BP B\times \BP C$. 
Now  consider  the \lq\lq limiting ideal\rq\rq. But how should one  take limits?  
If one works in the usual Hilbert scheme the $r$ points limit to some zero dimensional
scheme and one could take the span of that scheme. But for secant varieties one
is really taking the limit of the spans $\langle T_{1,\ep}\hd T_{r,\ep}\rangle$ in the Grassmannian
of $r$ planes and   the span of the limiting scheme can be strictly smaller than the limit
of the spans, so this idea does not work.

The answer is to work in the Haiman-Sturmfels multigraded Hilbert scheme  \cite{MR2073194}, 
which lives in a product of Grassmannians and keeps track of the entire Hilbert function rather than just the
Hilbert polynomial. The price one pays is that now one must allow unsaturated ideals. 
 
As with the border substitution method, one can insist that limiting ideal $I_0$ is Borel fixed,
which  for tensors with a large symmetry group  reduces 
in small multi-degrees
to a small search
 that has been exploited in practice.  

 Instead of single curve $E_\ep\subset G(r,A\ot B\ot C)$ limiting to a Borel
 fixed point, for each
 $(i,j,k)$ one gets a curve  in each
 $Gr(r,S^iA^*\ot S^jB^*\ot S^kC^*)$, each curve  limiting to a Borel fixed point {\it and}
 satisfying compatibility conditions. Here $Gr(r,V)$ is the Grassmannian of {\it codimension} $r$ subspaces in $V$.
 In this situation,  $E_\ep=(I_\ep)_{(111)}^\perp$.

The upshot is an algorithm that either produces all normalized candidate $I_0$'s
 or proves border rank $>r$. The caveat is that once one has a candidate $I_0$ one must determine
 if it actually came from a curve of ideals of $r$ distinct points.
 
 Using border apolarity, in \cite{CHLapolar} we proved numerous new matrix multiplication border rank lower bounds, 
 as well as determining the border rank of  the size three determinant polynomial considered as a tensor $\det_3\in 
 \BC^9\ot \BC^9\ot \BC^9$, whose importance for complexity is explained below.
 In particular, our results include the first nontrivial lower bounds for \lq\lq unbalanced matrix multiplication tensors\rq\rq ,
 something that was untouchable using previous methods.

\subsection{Border apolarity is subject to the cactus barrier}
In practice, one attempts to construct an ideal by building it up from low multi-degrees. The main restrictions
one obtains is when one has the ideal constructed up to  multi-degrees $(s-1,t,u)$, $(s,t-1,u)$ and $(s,t,u-1)$, 
and one wants to construct the ideal in degree $(s,t,u)$. In order that the construction may continue,
one needs that the natural symmetrization and  addition  map 
\be\label{stu}
I_{s-1,t,u}\ot A^*\oplus I_{s,t-1,u}\ot B^*\oplus I_{s,t,u-1}\ot C^*\ra S^sA^*\ot S^tB^*\ot S^uC^*
\ene
has image of codimension at least $r$. Here $I_{s-1,t,u}\subset S^{s-1}A^*\ot S^tB^*\ot S^uC^*$ denotes the component of the ideal
in multi-degree $(s-1,t,u)$ etc. (Here and in what follows, the direct sum is the abstract direct sum of vector spaces, so
there is no implied assertion that the spaces are disjoint when they live in the same ambient space.)

That is, the minors of the map of appropriate size must vanish. These are determinantal conditions
and therefore subject to the cactus barrier.

\begin{remark} Aside for the experts: J. Buczy{\'n}ski points out that not all components of the usual
Hilbert scheme contain ideals with generic Hilbert function. Thus in those situations, border apolarity may give 
better lower bounds on border rank than on cactus border rank.
\end{remark}
 
\subsection{Deformation theory}
Although  border apolarity alone cannot overcome the cactus barrier,   by placing calculations in the Haiman-Sturmfels multigraded
Hilbert scheme, it provides  a path to overcoming the cactus barrier. Namely one can apply the tools
of {\it deformation theory} (see, e.g., \cite{MR2583634} for an introduction) to determine if a candidate ideal is deformable to the ideal of a smooth
scheme.
Below, after motivating the problem,  I describe a first implementation of this in the   special case of tensors of minimal border rank.

\section{  Strassen's laser method for upper-bounding the exponent of matrix multiplication using
 tensors of minimal or near minimal  border rank}
\label{upperbnds}
The best way to prove an upper bound on matrix multiplication complexity 
would be to prove an upper bound for matrix multiplication directly. Fortunately for
algebraic geometry, Bini \cite{MR605920} showed that this is measured by the border rank 
of the matrix multiplication tensor. However, there has been little progress in this direction since
the early 1980's. Instead, border rank upper bounds for $\Mn$ have been proven using
indirect methods,  the most important papers are those of 
Sch\"onhage \cite{MR623057}, Strassen \cite{MR882307} and Coppersmith-Winograd
\cite{MR91i:68058}. The resulting technique is
called {\it Strassen's laser method}. Remarkably, it begins with a tensor of   minimal (or near minimal) border
rank, i.e., a concise tensor in $\BC^m\ot \BC^m\ot \BC^m$ of border rank $m$ or nearly $m$, and then
builds a large tensor from it, using its Kronecker powers defined below, and then,
using methods from probability and combinatorics, shows this large tensor admits a degeneration  to a large
matrix multiplication tensor.

For tensors $T\in A\ot B\ot C$ and $T'\in A'\ot B'\ot C'$, the {\it Kronecker product} of $T$ and $T'$ is the tensor $T\boxtimes T' := T \ot T' \in (A\ot A')\ot (B\ot B')\ot (C\ot C')$, regarded as $3$-way tensor. Given $T \in A \otimes B \otimes C$, the {\it Kronecker powers}  of $T$ are $T^{\boxtimes N} \in A^{\otimes N} \otimes B^{\otimes N} \otimes C^{\otimes N}$, defined iteratively. Rank and border rank are submultiplicative under Kronecker product: $\bfR(T \boxtimes T') \leq \bfR(T) \bfR(T')$, $\uR(T \boxtimes T') \leq \uR(T) \uR(T')$, and both inequalities may be strict.

Given $T,T' \in A \otimes B \otimes C$, we say that $T$ \emph{degenerates} to $T'$ if $T' \in \bar{GL(A) \times GL(B) \times GL(C) \cdot T}$, the closure of the orbit of $T$, the closures are the same
 in the Euclidean and   Zariski topologies.

Strassen's {\it  laser method} \cite{MR882307,MR664715} obtains upper bounds on $\om$ by
showing a certain explicit degeneration 
of a large Kronecker power of a    tensor $T$ satisfying certain combinatorial properties, 
admits a further  degeneration to a large matrix multiplication tensor.
Since border rank is non-increasing under degenerations and one has
an upper bound on $\ur(T^{\boxtimes N})$ inherited from the knowledge of $\ur(T)$, on obtains
an upper bound on the border rank of the large matrix multilplication tensor.

 An early  success of the laser method was with a tensor of border rank $m+1$,
now  called the {\it small Coppersmith-Winograd tensor} $T_{cw,q}\in (\BC^{q+1})^{\ot 3}$. Coppersmith and
Winograd showed that
for all $k$ and $q$,  \cite{MR91i:68058}
\begin{equation}\label{cwbndq}
\omega\leq  \log_q(\frac 4{27}(\uR(T_{cw,q}^{\boxtimes k}))^{\frac 3k}) .
\end{equation}
They used this when $q=8$ and the estimate  $\uR(T_{cw,q}^{\boxtimes k})\leq (q+2)^k$
to obtain $\omega\leq 2.41$. In particular, one could even potentially prove $\omega$ equaled two were
$\tlim_{k\ra\infty} (\uR(T_{cw,2}^{\boxtimes k}))^{\frac 3k} $ equal to $3$ instead of $4$.
Using an enhancement of border apolarity, with A. Conner and H. Huang, in \cite{CHLlaser} we  
 solved the longstanding
problem of determining $\uR(T_{cw,2}^{\boxtimes 2})$. 
Unfortunately for matrix multiplication upper bounds, we proved that $\uR(T_{cw,2}^{\boxtimes 2})=4^2$.
Previously, just using Koszul flattenings, analogous (and even higher Kronecker power)
 results for other small Coppersmith-Wingorad tensors
 were obtained  with A. Conner, F. Gesmundo, and E. Ventura \cite{MR4356227}.

A more intriguing tensor is the \lq\lq skew-cousin\rq\rq  of the small Coppersmith-Winograd tensor 
$T_{skewcw,q}$ occuring in odd dimensions, which similarly satisfies
for all $k$ and even $q$,  \cite{MR4356227}  
\begin{equation}\label{skewcwbndq}
\omega\leq  \log_q(\frac{4}{27}(\uR(T_{skewcw,q}^{\boxtimes k}))^{\frac{3}{k}}) .
\end{equation}
Again, the $q=2$ case could potentially be used to prove the exponent is two.
Here one begins with a handicap, as    $\ur(T_{skewcw,2})=5>4$,  but with
A. Conner and A. Harper, using border apolarity for the lower bound and numerical search
methods for the upper bound, in  \cite{CHLapolar} we showed $\ur(T_{skewcw,2})^{\boxtimes 2}=17<25$.
Unfortunately $17>16$. The  problem of determining the border rank of the cube
remains.

It is worth remarking that $T_{cw,2}^{\boxtimes 2}$ is isomorphic to the size three permanent polynomial
considered as a tensor and $T_{skewcw,2}^{\boxtimes 2}$ is 
 isomorphic to the size three determinant polynomial \cite{MR4356227}.

 The best bounds on the exponent were obtained using the laser method applied to the
 {\it big Coppersmith-Winograd tensor}  $T_{CW,q}$, which has minimal border rank.  
 However, 
 barriers to future progress using the laser method applied to this tensor have been discovered, first in \cite{MR3388238},
 and then in numerous follow-up works. In particular, one cannot prove $
 \om<2.3$ using 
 $T_{CW,q}$ in the laser method.
 A geometric interpretation of the barriers is given in \cite{MR3984631}. 
 
 Very recently, at an April 2022 workshop on geometry and complexity theory in Toulouse,
 France,  J. Jelisiejew and M. Micha{\l}ek announced a path to improving the laser method. Their
 observation was that the border rank estimate for the \lq\lq certain degeneration\rq\rq\ of $T^{\boxtimes N}$ 
 in the laser method mentioned above can be improved! The proof exploits properties of the algebra
 associated to  $T^{\boxtimes N}$ (discussed in \S\ref{bindingsect} below) that persist under the degeneration.
 
Even without that recent developement, other minimal border rank tensors could potentially prove $\om<2.3$
with the standard laser method. In fact in \cite[Cor 4.3]{hoyois2021hermitian} and 
 \cite[Cor 7.5]{2019arXiv190909518C} it was observed that among
 tensors that are $1_A$, $1_B$ and $1_C$ generic (such are called {\it  $1$-generic tensors}), $T_{CW,q}$ is
 the \lq\lq worst\rq\rq \ for the laser method in the sense that any bound one can prove using $T_{CW,q}$
 can also be proved using any other minimal border rank $1$-generic tensor. This provides strong motivation
 to understand tensors of minimal border rank. A second motivation is that it can serve as a  case study
 for the implementation of deformation theory to overcome the cactus barrier.

\section{Classical and neo-classical equations for tensors of  minimal border rank}
 
 Before discussing recent developments for tensors of minimal border rank, I explain the previous
 state of the art.
I already have discussed Strassen's equations and Koszul flattenings. What follows are additional conditions.

\subsection{The equations of \cite{LMsec,LMsecb}} Several modules of equations were found in \cite{LMsec,LMsecb}
using representation theory and variants of Strassen's equations. Many of these still lack a geometric
interpretation.

\subsection{The flag condition}
If
$\ur(T)=m$   there exists a flag $A_1\subset \cdots \subset A_{m-1}\subset A$ such that
for all $j$,  $T(A_j^*)\subset \s_j(Seg(\BP B\times \BP C))$.
This has been observed several times, dating back at least to \cite{BCS}. 
Note that to convert this condition to polynomial conditions, one would have to use elimination theory, even
for the first step that there exists a line $A_1^*$ such that  $\BP T(A_1^*)\in  Seg(\BP B\times \BP C)$.
The flag condition was essential to the results in \cite{CHLlaser}.

\subsection{The End-closed condition} \label{endclosedpoly}
Gerstenhaber \cite{MR0132079} observed the following: Let $\langle x_1\hd x_m\rangle\subset \tend(\BC^m)$ be a limit of 
  spaces of simultaneously diagonalizable matrices.   
Then $\forall i,j$, $x_ix_j\in \langle x_1\hd x_m\rangle$.   
Call this the 
\lq\lq End closed condition\rq\rq . To express the condition  as polynomials,
let $\{\a_i\}$ be a basis of $A^*$, with $\a_1$ chosen to maximize the rank of $T(\a_1)$, then 
  for all $\a',\a''\in A^*$,  the End-closed condition is 
\be\label{bigenda1gen}
(T(\a')T(\a_1)^{\ww m-1}T(\a'') )
\ww T(\a_1) \ww \cdots \ww  T(\a_m) 
=0\in \La{m+1}(\tend(C)).
\ene
These are 
  polynomials of degree $2m+1$. 
When $T$ is $1_A$-generic and one takes $\a_1$ such that $\trank(T(\a_1))=m$,
 these correspond to $T(A^*)T(\a_1)\inv\subset \tend(C)$ being closed
under composition of endomorphisms.

\subsection{The symmetry Lie algebra condition}
Let $\fg=\fgl(A)\op \fgl(B)\op \fgl(C)$. Let $\hat \fg_T=\{ X\in \fg\mid X.T=0\}$. (This is the pullback of
the symmetry Lie algebra of $T$ to $\fgl(A)\op \fgl(B)\op \fgl(C)$.)
With $T$ understood, write
$\fg_{AB}=\{ X\in \fgl(A)\op \fgl(B)\mid X.T=0\}$ and similarly for $\fg_{BC},\fg_{AC}$. 

A concise tensor of rank $m$, $\Mone^{\op m}$,  has $\tdim \hat \fg_{\Mone^{\op m}}=2m$ and
$\tdim \fg_{AB}=\tdim \fg_{AC}=\tdim \fg_{BC}=m$.
The dimension of the symmetry Lie-algebra is semi-continuous under degenerations, thus 
if $T$ is of minimal border rank $\tdim \hat \fg_{T}\geq 2m$ and $\tdim \hat \fg_{AB}\geq m$ and permuted statements.

Computing these dimensions amounts to determining the dimension of the kernel of a linear map.
Precisely to check if  $\tdim \hat \fg_{T}\geq 2m$ are equations of degree $3m^2-2m+1$
and $\tdim \fg_{AB}\geq m$ are equations of degree $2m^2-m+1$.

\section{The 111-equations and first consequences}\label{111test}

\subsection{The $(111)$-equations}

The {\it 111-equations} are the rank conditions on the map \eqref{stu} when $(s,t,u)=(1,1,1)$
and one is testing for border rank $m$. Note that in this case there are no choices
for the ideal in degrees $(110),(101),(011)$,  so they are really polynomial equations. 
 These equations first appeared  in \cite[Thm 1.3]{MR4332674}.

The  $111$-equations for concise tensors of minimal border rank  may be rephrased as the requirement that 
\be\label{triplei}\tdim((T(A^*)\ot A)\cap (T(B^*)\ot B) \cap (T(C^*)\ot C))\geq m.
\ene

A special case of the $111$-equations are the two-factor $111$-equations, which
have a natural geometric interpretation and  are easier to implement
because a pairwise intersection 
can be computed using inclusion-exclustion:
Given
subspaces $X_1,X_2,X_3$ of a vector space $V$, by inclusion-exclusion
$\tdim (X_i\cap X_j)=\tdim(X_i)+\tdim(X_j)-\tdim \langle X_i,X_j\rangle$.

Thus the two-factor $(111)$-test may be computed by checking if 
$\tdim \langle T(A^*)\ot A^*, T(B^*)\ot B^*\rangle\geq 2m^2-m+1$ and permuted statements.
These are equations of degree $2m^2-m+1$ in the $T^{ijk}$.  Notice that if $(X,Y)\in \fg_{AB}$, i.e., $X.T=-Y.T$, then
$(X,-Y)$ gives rise to an element of $(T(A^*)\ot A)\cap (T(B^*)\ot B)$, i.e., the two factor 111-tests are equivalent to
the dimension requirements on $\fg_{AB},\fg_{AC},\fg_{BC}$ for minimal border rank.

More generally, the full     $(111)$-equations   may also be understood as a generalization of
the lower bound on $\tdim(\hat\fg_T)$, where one  
not just bounds dimension, but restricts the structure of $\hat\fg_T$ as well.

To compare the 111-equations with other previously known equations, we have:

\begin{proposition}\cite[Prop. 1.1, Prop. 1.2]{https://doi.org/10.48550/arxiv.2205.05713} The 111-equations imply both Strassen's equations and the End-closed equations.
The 111-equations do not always imply the $p=1$ Koszul flattening equations.
\end{proposition}

Consider  the situation of  a concise tensor where each of the associated spaces
of homomorphisms is of bounded rank $m-1$.   Strassen's equations do allow some
assertions in this situation. A normal form for such tensors was proved
by S. Friedland \cite{MR2996364}.  This normal form was generalized in
\cite[Prop. 3.3]{https://doi.org/10.48550/arxiv.2205.05713}
by using the 111-equations
instead of Strassen's equations. (In fact this generalized normal form allowed the proof
that the 111-equations imply Strassen's equations and the End-closed equations.) These normal forms respectively allowed the
characterization of concise tensors of minimal border rank when $m=4$ and $m=5$, in fact
S. Friedland was even able to resolve the non-concise $m=4$ case using additional equations
he developed, solving the set-theoretic \lq\lq salmon prize problem\rq\rq\  posed by E. Allman.

Recall that Strassen's equations and the End-closed equations are trivial when a tensor gives
rise to three linear spaces of bounded rank at most $m-2$. 
The 111-equations do not share this defect. We are currently implementing them for
such spaces of tensors.
(The $p=1$ Koszul flattenings
are not trivial in this setting, we have yet to determine their utility for bounded rank $m-2$ situations.)

\section{Deformation theory for tensors of minimal border rank}

For tensors $T\in \BC^m\ot \BC^m\ot \BC^m$ satisfying genericity conditions, one has natural algebraic structures associated 
to them that can  be utilized to help determine if they have minimal border rank. 

\subsection{Binding tensors and algebras}\label{bindingsect}
Say $T\in A\ot B\ot C$ is $1_A$ and $1_B$ generic with $T(\a_1): B^*\ra C$ and $T(\b_1):A^*\ra C$ of full rank.
(A tensor that is at least two of $1_A$, $1_B$ or $1_C$ generic is called {\it binding}.)
Use their inverses to obtain a tensor isomorphic to $T$, which I abuse notation and also denote by $T$,
$T\in C^*\ot C^*\ot C$, i.e.,  a bilinear map $T: C\times C\ra C$,  which gives $C$ the structure of an algebra with left identity $\a_1$ and right
identity $\b_1$.

If $T$ satisfies the  $A$-Strassen  equations then it is isomorphic to a partially symmetric tensor, see Proposition
\ref{prop1},  and the associated algebra is abelian. Conversely, given such an algebra, one obtains
its structure tensor.

Explicitly, 
let $\ci\subset \BC[x_1\hd x_n]$ be an ideal whose zero set in affine space is finite,
more precisely so that $\cA_\ci:=\BC[x_1\hd x_n]/\ci$ is a finite dimensional algebra of dimension $m$.
(This will be the case, e.g.,  if the zero set consists of $m$ distinct points each counted with multiplicity one.)
Let $\{p_I\}$ be a basis of $\cA_\ci$ with dual basis $\{p_I^*\}$
We can write the structural tensor of $\cA_\ci $ as
$$
T_{\cA_\ci}=\sum_{p_I,p_J\in \cA_\ci} p_I^*\ot p_J^*\ot (p_Ip_J \tmod \ci).
$$

Then \cite{MR3578455}   shows that a binding tensor $T$ that is the structure tensor of a smoothable algebra
is of minimal border rank, 
i.e., the tensor $\Mone^{\op m}$
degenerates  to  $T$, where
 $\Mone^{\op m}$ is the tensor 
 whose associated algebra $\cA_{\Mone^{\op m}}$ comes from the ideal of $m$ distinct points.
The key step is showing that in this situation $T\in \ol{GL(A)\times GL(B)\times GL(C)\Mone^{\op m}}$
if and only if (using the above identifications) $T\in \ol{  GL(C)\Mone^{\op m}}$.

Thus one may utilize deformation theory on the Hilbert scheme of points to determine if
a binding tensor satisfying  the $A$-Strassen equations has minimal border rank.
In particular, such tensors automatically are of minimal border rank when $m\leq 7$ \cite{MR2579394}.

\subsection{$1$-Generic tensors: Gorenstein algebras}
Now say $T$ is $1_A$, $1_B$, and $1_C$ generic (such tensors are called {\it  $1$-generic}) and satisfies the $A$-Strassen  equations. We have $\g_1\in C^*$
such that $T(\g_1)\in \tend(C)$ is invertible. What extra structure do we obtain?

Recall that an algebra $\cA$ is {\it Gorenstien} if there exists   $f\in \cA^*$ such that
any of the following equivalent conditions holds:

\noindent 1) $T_{\cA}(f)\in \cA^* \otimes \cA^*$ is of full rank,

\noindent 2) the pairing $\cA\otimes \cA\ra \BC$  given by $(a, b)\mapsto f(ab)$ is non-degenerate,

\noindent 3) $\cA f = \cA^*$.
 
 
Thus $f=\g_1$ above  tells us $\cA_T$ is Gorenstein by (1). By the second assertion in Proposition \ref{prop1},
$T$ is moreover symmetric.

In particular, such $T$ is of minimal border rank when $m\leq 13$ \cite{MR3404648}.
For an algorithm that resolves the $m=14$ case,  see \cite{https://doi.org/10.48550/arxiv.2007.16203}.

The additional algebraic structure of being Gorenstein   makes the deformation theory
easier to implement.

\begin{example}  Consider $\cA=\BC[x]/(x^2)$, with basis $1,x$, so
$$
T_{\cA}=1^*\ot 1^*\ot 1 + x^*\ot 1^*\ot x + 1^*\ot x^*\ot x.
$$
Writing $e_0=1^*$, $e_1=x^*$ in the first two factors and
$e_0=x$, $e_1=1$ in the third,  
$$
T_{\cA}=e_0\ot e_0\ot e_1+e_1\ot e_0\ot e_0+e_0\ot e_1\ot e_0
$$
That is, $T_{\cA}=T_{WState}$ is a general tangent vector to $Seg(\BP A\times \BP B\times \BP C)$.

 \end{example}
  
\begin{example}[The big Coppersmith-Winograd tensor]  
Consider the algebra
$$
\cA_{CW,q}=\BC[x_1\hd x_q]/(x_ix_j, x_i^2-x_j^2, x_i^3, \ i\neq j)
$$
Let $\{ 1, x_i, [x_1^2]\}$ be a basis of $\cA$, where $[x_1^2]=[x_j^2]$ for all $j$. 
Then
\begin{align*}T_{\cA_{CW,q}}
=&
1^*\ot 1^*\ot 1+ \sum_{i=1}^q (1^*\ot x_i^* \ot x_i + x_i^*\ot 1^*\ot x_i )\\
&
+ x_i ^*\ot  x_i^* \ot [x_1^2] + 1^*\ot [x_1^2]^*\ot [x_1^2] + 
[x_1^2]^*\ot 1^*\ot [x_1^2].
\end{align*}
Set $e_0=1^*$, $e_i=x_i^*$, $e_{q+1}=[x_1^2]^*$
in the first two factors and
$e_0=[x_1^2]$, $e_i=x_i$, $e_{q+1}=1$ in the third to obtain
\begin{align*}
T_{\cA_{CW,q}}=&T_{ CW,q}
=e_0\ot e_0\ot e_{q+1}+ \sum_{i=1}^q (e_0\ot e_i\ot e_i+
e_i\ot e_0\ot e_i + e_i\ot e_i\ot e_0)\\
&+
e_0\ot e_{q+1}\ot e_0+e_{q+1}\ot e_0\ot e_0, 
\end{align*}
which is the usual expression for  the big Coppersmith-Winograd tensor. 
\end{example}

\subsection{$1_*$-generic tensors: modules and the ADHM correspondence}
 When $\tdim A=\tdim B=\tdim C=m$, one says $T$ is {\it $1_*$-generic} if it is 
$1_A$ or $1_B$ or $1_C$ generic. 

Consider the case of  a tensor that is   $1_A$-generic but not binding. What structure can we associate to it? Fixing $\a_1$ as
above we obtain $T\in A\ot C^*\ot C$, i.e., $T(A^*)T(\a_1)\inv\subset \tend(C)$, and if Strassen's equations
are satisfied, we have an abelian subspace, and if furthermore the End-closed condition holds,
we may think of this space as defining an algebra action on $\tend(C)$, which we may lift to
an action of the polynomial ring $S:=\BC[y_2\hd y_m]$ by $y_s(c):=T(\a_s)T(\a_1)\inv c$.
(The choice of indices is deliberate, as $T(\a_1)T(\a_1)\inv=\Id_C$ corresponds to $1\in S$.)

That is, the vector space $C$ becomes a module over the polynomial ring. This association
is called the {\it ADHM correspondence} in   \cite{jelisiejew2021components},  after  \cite{MR598562}. This leads one
to deformation theory in the Quot scheme that parametrizes such modules.

This correspondence allowed Jelisiejew, Pal and myself \cite{https://doi.org/10.48550/arxiv.2205.05713} to characterize concise  $1_*$-generic tensors
of border rank $\leq 6$ as the zero set of Strassen's equations and the End-closed equations, and also
as the zero set of the $111$-equations. Strassen's equations, the $111$-equations and the End-closed 
equations fail to characterize minimal border rank tensors  when $m\geq 7$.

\subsection{Concise tensors: the 111-algebra and its modules}
Now say we just have a concise tensor. Previously there had not been any algebraic structure
available for studying such tensors. Moreover, as remarked above, both Strassen's equations and
the End-closed equations are trivially satisfied for such tensors when the three associated spaces
of homomorphisms are of rank bounded above by $m-2$.  Despite this, the 111-equations
still give strong restrictions in these cases. I now explain that they also allow the implementation
of deformation theory even in this situation. 

 For   $\Amat\in \tend(A) = A^*\ot A$, let
    $\Amat\acta T$  
    denote the corresponding element of $T(A^*)\ot A$. Explicitly, if $\Amat =
    \alpha\ot a$, then $\Amat \acta T := T(\alpha)\ot a$ and the map $(-)\acta
    T\colon
    \tend(A)\to A\ot B\ot C$ is extended linearly.
  
    \begin{definition}\cite[Def. 1.9]{https://doi.org/10.48550/arxiv.2205.05713}
        Let $T$ be a concise tensor.
        We say that a triple $(\Amat, \Bmat, \Cmat)\in \tend(A)
        \times\tend(B)\times \tend(C)$ \emph{is compatible with} $T$ if
        $\Amat\acta T = \Bmat \actb T = \Cmat \actc T$.
        The \emph{111-algebra} of $T$ is the set of triples
        compatible with $T$.  We denote this set by $\alg{T}$.
    \end{definition}
    
    Thus a compatible triple gives a point in the triple intersection \eqref{triplei}.
    The name 111-{\it algebra}  is justified by the following theorem:
    
      \begin{theorem}\label{ref:111algebra:thm}\cite[Thm. 1.10]{https://doi.org/10.48550/arxiv.2205.05713}
        The 111-algebra of a concise tensor $T\in A\ot B\ot C$ is a
        commutative unital
        subalgebra of $\tend(A)\times \tend(B) \times \tend(C)$ and its projection to any factor is injective.
    \end{theorem}

Using the 111-algebra, one obtains  four consecutive obstructions for a concise tensor to be of minimal border
rank \cite{https://doi.org/10.48550/arxiv.2205.05713}:
\begin{enumerate}
    \item\label{it:abundance} $\tdim(  \alg{T})\geq m$. For what the next three conditions, assume equality holds.
    \item\label{it:cactus}   $\alg{T}$ must  be smoothable.   
    \item\label{it:modulesPrincipal}  Using the 111-algebra, $A,B,C$ become
    modules for it and the polynomial ring $S$.  These three  $S$-modules, $\ul A,\ul B,\ul C$ (where
    the underline is there to emphasize  their module structure)   must  lie in the principal component
    of the Quot scheme.
    \item\label{it:mapLimit} there exists a surjective module homomorphism $\ul A\ot_{\alg{T}} \ul B\to \ul C$ associated to $T$  and this homomorphism must  be a  limit of
       module homomorphisms  $\ul A_\ep\ot_{\cA_\ep} \ul B_\ep \to \ul C_\ep$
       for a
        choice of smooth algebras $\cA_\ep$ and semisimple modules $\ul A_{\ep}$, $\ul B_{\ep}$, $\ul C_{\ep}$. 
\end{enumerate}

 \section{New proofs of existing results}  \label{endclosedpf}
 In this section I present  two significantly simpler proofs than the original that binding tensors satisfying Strassen's equations
 are partially symmetric and the original, more elementary proof that binding tensors satisfying Strassen's equations
 automatically satisfy the End-closed condition. These   proofs were obtained in conversations with
 J. Jelisiejew and M. Micha{\l}ek.

Let $A,B,C\simeq \BC^m$ and let $T\in A\ot B\ot C$ be $1_A$-generic. Say $\trank (T_A(\a_0))=m$.

Note the tautological identities:
$T(\a,\b)=T_A(\a)\b=T_B(\b)\a$.

The $A$-Strassen equations  for minimal border rank say that for all $\a_1,\a_2\in A$,
$$T_A(\a_1)T_A(\a_0)\inv T_A (\a_2)= T_A(\a_2)T_A(\a_0)\inv T_A (\a_1).
$$

\begin{proposition}\label{prop1}\cite{MR3682743}  Let $T$ be $1_A$ and $1_B$-generic and satisfy the $A$-Strassen equations.
Then $T$ is isomorphic to a tensor in $S^2C^*\ot C$. If $T$ is $1$-generic then it is isomorphic to a 
symmetric tensor.
\end{proposition}

Here are two proofs:

\begin{proof} Assume $T(\a_0)\in B\ot C$ and $T(\b_0)\in A\ot C$ are of full rank. 
Define $\tilde T\in C^*\ot C^*\ot C$ by $\tilde T(c_1,c_2):=T(T_B(\b_0)\inv c_1, T_A(\a_0)\inv c_2)$.

Set $\a_1=T_B(\b_0)\inv c_1$, $\a_2=T_B(\b_0)\inv c_2$ so
\begin{align*}
\tilde T(c_1,c_2)&=T(\a_1,T_A(\a_0)\inv T_B(\b_0)\a_2)\ \ {\rm definition}\\
&=T(\a_1,T_A(\a_0)\inv T_A(\a_2)\b_0)\ \ {\rm taut. id.}\\
&=T_A(\a_1)T_A(\a_0)\inv T_A(\a_2)\b_0\ \ {\rm taut. id.}\\
&=T_A(\a_2)T_A(\a_0)\inv T_A(\a_1)\b_0\ \ {\rm Strassen}\\
&=\tilde T(c_2,c_1)\ \ {\rm taut. id.}
\end{align*}
The second assertion follows as $\FS_3$ is generated by the transpositions $(1,2)$ and $(1,3)$.
\end{proof}

\begin{proof} Under the hypotheses
$T_A^{\ww 1}: A \ot B^*\ra \La 2 A\ot C$ has rank $m^2-m$ and the $1_B$-genericity
condition assures that the $m$-dimensional kernel contains an element of full rank,
$\psi: A\ra B$, which makes $(\psi\ot \Id_{B\ot C})(T)\in S^2B^*\ot C$. The second assertion follows similarly.
\end{proof}

Note that Proposition \ref{prop1} implies the $B$-Strassen equations are satisfied as well.

The following proposition appeared in \cite{https://doi.org/10.48550/arxiv.2205.05713} with a less elementary proof. Below   is the original
proof. 

\begin{proposition}  If $T$ is $1_A$ and $1_B$ generic  and satisfies the $A$-Strassen equations, then $T(A^*)T(\a_0)\inv
\subset \tend(C)$ satisfies the End-closed condition.
\end{proposition}
\begin{proof}
We need to show that for all $\a_1,\a_2$, that, there exists $\a'$ such that
$$T_A(\a_1)T_A(\a_0)\inv T_A(\a_2)T_A(\a_0)\inv =T_A(\a')T_A(\a_0)\inv.
$$

It is sufficient to work with $\tilde T\in S^2C^*\ot C$. Here,   by symmetry $\tilde T_A(c)=\tilde T_B(c)=: \tilde T_{C^*}(c)$.
We claim
$\tilde T_{C^*}(c_1)\tilde T_{C^*}(c_2)=\tilde T_{C^*}(\tilde T(c_1,c_2))$.

To see this
\begin{align*}\tilde T_{C^*}(c_1)\tilde T_{C^*}(c_2)(c)&=\tilde T_{C^*}(c_1)\tilde T (c_2,c)\ \ {\rm taut.}\\
&=  T_{C^*}(c_1)\tilde T (c ,c_2)) \ \  {\rm sym.}\\
&=\tilde T_{C^*}(c_1) \tilde T_{C^*} (c  )(c_2)\ \ {\rm  taut.} \\
&=\tilde T_{C^*}(c) \tilde T_{C^*} (c_1 )(c_2)\ \ {\rm Strassen} \\
&=\tilde T(c, \tilde T(c_1,c_2))\ \ {\rm   taut.}\\
&=\tilde T( \tilde T(c_1,c_2),c)\ \  {\rm sym.}\\
&=\tilde T_{C^*}(\tilde T(c_1,c_2))(c) \ \ {\rm   taut.}
\end{align*}

If $T$ is binding with $T_A(\a_0),T_B(\b_0)$ of full rank, we may take $T(\a_0,\b_0)\in C$
as the generator of the algebra.
We need to show that in this situation
$$
T_A(\a_0)T_A(\a_0)\inv T(\a_0,\b_0), T_A(\a_1)T_A(\a_0)\inv T(\a_0,\b_0)\hd T_A(\a_n)T_A(\a_0)\inv T(\a_0,\b_0)
$$
is a basis of $C$.
But this is
$$
\{T_A(\a_0) \b_0, T_A(\a_1) \b_0\hd T_A(\a_n) \b_0\}=
\{T_B(\b_0)\a_0\hd T_B(\b_0)\a_n\}
$$
Since $T_B(\b_0): A^*\ra C$ is of maximal rank, the basis $\a_0\hd \a_n$ of $A^*$ will map to 
a basis of $C$.
\end{proof}

  \bibliographystyle{amsplain}

\bibliography{Lmatrix}

\end{document}